\documentclass[a4paper, onecolumn]{IEEEtran}
\usepackage{amssymb}
\usepackage{amsmath}	
\usepackage{mathtools}
\usepackage{amsthm}
\usepackage{tikz,pgf,pgfplots,circuitikz}
\usepackage{amsmath}

\usetikzlibrary{shapes.geometric}  
\usetikzlibrary{patterns}
\usetikzlibrary{decorations.text}
\usetikzlibrary{%
    decorations.pathreplacing,%
    decorations.pathmorphing%
}
\usetikzlibrary{mindmap,trees}
\usetikzlibrary{calc}
\usetikzlibrary{fadings}
\usetikzlibrary{decorations.pathreplacing}
\newcommand{\colora}{green!40!blue}   
\newcommand{\colorb}{green!60!blue}

\newcommand{\colore}{blue!60}

\usepackage{listings}
\title {Stability Analysis in Type-A Wind Turbines: A Tutorial}

\author{\IEEEauthorblockN{Alejandro Garc\'es}\\
\IEEEauthorblockA{Department of Electric Power Engineering\\ Universidad Tecnol\'ogica de Pereira (UTP)\\
AA: 97 - Post code: 660003 \\
Pereira, Colombia}}

\begin{document}
\maketitle
\begin{abstract}
	This tutorial presents a simple yet accurate transient stability analysis for Type-A wind turbines. The paper is presented in a tutorial way and therefore it includes some scripts in Matlab/muPad which demonstrate the main concepts. 
	
\end{abstract}
\section{Introduction}

Early wind turbines used to be integrated to grid thought a synchronous machine. In that type of machine, rotational mechanical speed is proportional to the frequency of the stator flux.   This frequency is imposed by the system and hence, the mechanical speed must be constant even if the wind velocity changes.  However, the maximum efficiency of the turbine is achieved if the rotational speed is proportional to the wind velocity with a well defined proportional constant.  A synchronous machine achieves maximum efficiency only at nominal value of the wind velocity, an operation point that is achieved very few times.  As consequence, the synchronous machine is not the best alternative for integrating wind turbines.

\begin{figure}
\centering
\footnotesize
\begin{tikzpicture}[x=0.30mm,y=0.30mm]	
	\draw[gray,left color=gray!20, right color=gray]
  (10,13) to[out=180, in=60] +(-25,-13) to[out=-60, in=180] +(25,-13) -- cycle;
  
   \draw[gray,left color=gray!20, right color=gray] (0,10) -- +(-3,20) -- +(0,100) to[out=60, in=90] +(2,100) --
 +(6,-70) to[out=-80, in=90] +(2,-100) -- cycle;
   \draw[gray,left color=gray!20, right color=gray] (0,-10) -- +(-3,-20) -- +(0,-100) to[out=-60, in=-90] +(2,-100) --
 +(6,70) to[out=80, in=-90] +(2,100) -- cycle;

 \draw [gray,fill=gray] (10,-2) rectangle (80,4);
% gear box
	\draw[-, gray,thick,fill = white] (20,-10) rectangle  +(30,20);
	\draw[-,gray,fill] (35,10) arc (360:270:15) -- +(0,15) -- cycle;
	\foreach \a in {0,-10,-20,-30,-40,-50,-60,-70,-80} \draw[gray,fill] (20,10) -- +(\a:18)  -- +(\a-5:18) -- cycle;
	\foreach \a in {100,110,120,130,140,150,160,170,180} \draw[gray,fill] (50,-10) -- +(\a:18)  -- +(\a-5:18) -- cycle;
	\draw[-,gray,fill] (35,-10) arc (180:90:15) -- +(0,-15) -- cycle;
% generador
 	\draw[-, \colora,thick,fill = white] (70,0) circle (12);
	\draw[-, \colora,thick,fill = \colora!30] (70,0) circle (10);
	\node at (70,0) {G}; 	
% trafo		
 	\draw[-,\colora,thick] (150,0) circle (10);
 	\draw[-,\colora,thick] (160,0) circle (10);
 	\draw[-,\colora,thick] (82,0) -- (140,0);
 	\draw[-,\colora,thick] (170,0) -- (190,0);
 	\draw[-,\colora,very thick] (190,-10) -- +(0,20);
% Arrancador
	\draw[-,\colora, thick,fill=white] (90,-5) rectangle +(20,10);
% capacitor
 	\draw[-,\colorb,thick] (120,0) -- +(0,-20);
 	\draw[-,\colorb,thick] (115,-20) -- +(10,0);
 	\draw[-,\colorb,thick] (115,-23) -- +(10,0);
 	\draw[-,\colorb,thick] (120,-23) -- +(0,-20);
       \draw[-,\colorb,thick] (117,-43) -- +(6,0); 	 	
       \draw[-,\colorb,thick] (118,-44.5) -- +(4,0); 	 	
       \draw[-,\colorb,thick] (119,-46) -- +(2,0); 	 	
% Arrancador 	

% Texto	
	\draw[-latex] (60,-80) node[right] () {Induction machine} to[out=180,in=-90] (70,-5) ;	
	\draw[-latex] (120,80) node[right] () {Gearbox} to[out=180,in=90] (40,10) ;
	\draw[-latex] (120,40) node[right] () {starter} to[out=180,in=90] (100,5) ;
	\draw[-latex,\colorb] (130,-60) node[right] () {Compensation} to[out=180,in=0] (125,-22) ;
	
\end{tikzpicture}
\caption{Schematic diagram of a Type-A wind turbine}
\label{fig:turbinatipoA}
\end{figure}
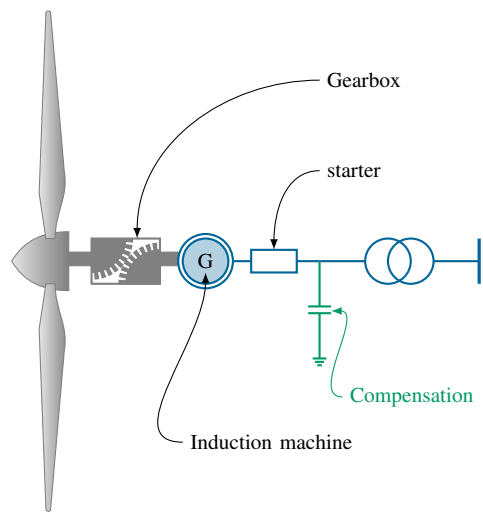

A second generation of wind turbines was developed in the 70's which was based in the asynchronous machine, usually an squirrel carroll induction generator as shown in Fig \ref{fig:turbinatipoA}.  This type of configuration is usually called \textit{Type-A wind turbine} \cite{libroolimpo}.  Its main advantage is the independence between the rotational speed and the frequency of the grid. This allows a better exploitation of the available power of the wind.  However, the use of asyncronous generators directly connected to the power grid entails some challenges related to the stability of the system.

Induction machines consume reactive power even if operated as generator.  This is because the machine requires to be magnetized before start generating power.   Consequently, a capacitive compensation is required in terminals of the machine.  Active compensation (i.g an STATCOM) is also a possibility though much more expensive.  In addition, a system which limits the current during starting up is required.  

Type-A wind turbines represents only $15\%$ of the wind parks installed around the world.  Its use is usually limited to medium and low power.  Therefore, it is an important element in the micro-grid concept.
 
\section{System modeling}

\subsection{Asynchronous generator}

Type-A wind turbines are integrated to the system by using induction machines which have a completely different behavior than  synchronous generators since their magnetization depend on the induction phenomena.  This means, the field in the rotor is never synchronized with the field in the stator in order to obtain an electrical torque. The difference between the rotation speed of these two fields is named slip ($s$) which can be defined as follows:

\begin{equation}
	s = \frac{w_{r}-w_{s}}{w_{s}}
\end{equation}

where $w_{s}$ is the angular electrical synchronous frequency and $w_{r}$ is the rotor electrical angular frequency.  Notice that under this definition $s>0$ correspond to the generator operation since the rotor field leads the stator field.    In per unit we can say that $s = w_{r}-1$ and constitutes the state variable for the model of the machine.

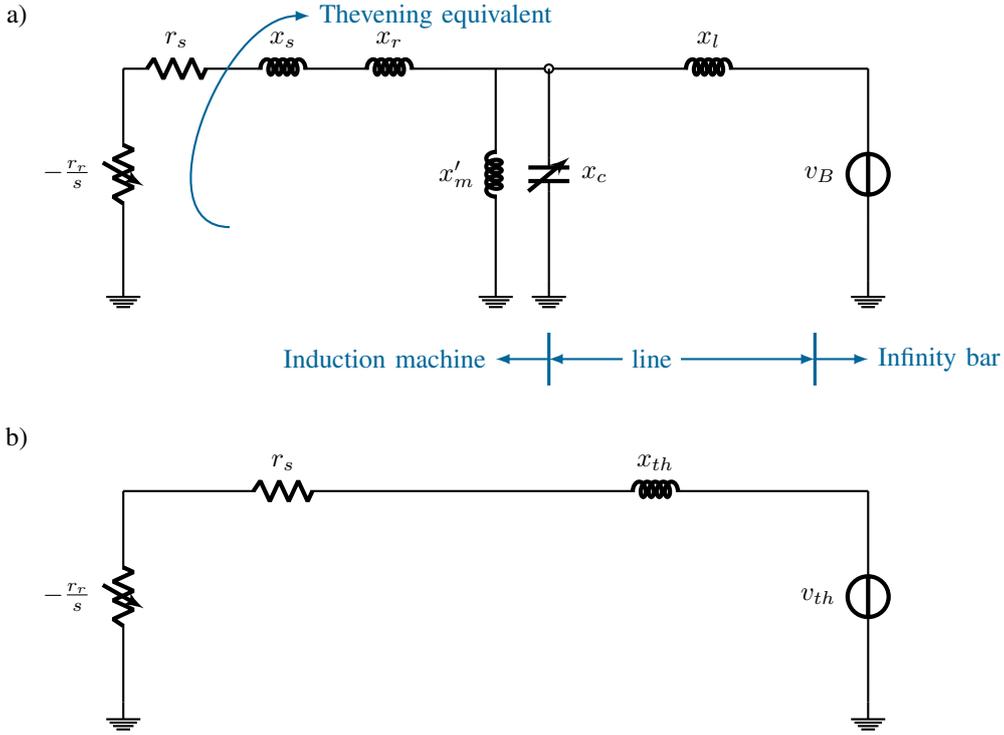
\begin{figure}
\centering 
\begin{tikzpicture}[scale=0.7]
\ctikzset{bipoles/length=0.90cm}
\draw[color=black, thick] (0,0) node[ground]{} 
                                                  to [vR, l = $-\frac{r_{r}}{s}$] (0,4)  
                                                  to [R, l = $r_{s}$] (2,4) 
                                                  to [L, l = $x_{s}$] (4,4)
                                                  to [L, l = $x_{r}$] (6,4)
                                                  to [short, -o] (8,4) 
                                                  to [L, l=$x_{l}$] (14,4);
\draw[color=black, thick] (14,0)node[ground]{} to [V, l=$v_{B}$] (14,4);                    
\draw[color=black, thick] (8,0) node[ground]{} 
                                                  to [vC, l_ = $x_{c}$] +(0,4)  ;                                                                                
\draw[color=black, thick] (7,0) node[ground]{} 
                                                  to [L, l = $x_{m}'$] (7,4);      
\draw[color=\colora, -latex, thick] (2,1) to[out=180, in=180] (3.5,5) node[right] {Thevening equivalent};                                                   
\draw[color=\colora, -latex, thick]  (8,-1.5) -- +(-1,0) node[left] {Induction machine};                                                
\draw[color=\colora,  very thick] (8,-2) -- +(0,1);
\draw[color=\colora,  very thick] (13,-2) -- +(0,1);
\draw[color=\colora, latex-latex, thick]  (8,-1.5) -- (10.5,-1.5) node[left, fill=white] {line} -- (13,-1.5);                              
\draw[color=\colora, -latex, thick]  (13,-1.5) -- +(1,0) node[right] {Infinity bar};                                                          

\node at (-2,5) {a)};
\node at (-2,-3) {b)};

\draw[color=black, thick] (0,-8) node[ground]{} 
                                                  to [vR, l = $-\frac{r_{r}}{s}$] +(0,4)  
                                                  to [R, l = $r_{s}$] +(6,0) 
                                                  to [L, l = $x_{th}$] +(8,0);
\draw[color=black, thick] (14,-8) node[ground]{} 
                                                  to [V, l = $v_{th}$] +(0,4);
\end{tikzpicture}
\caption{Induction machine equivalent circuit a). Complety model including external system b) thevening equvalent}
\label{fig:circuito_maquina_induccion}
\end{figure}

A simplified model of the asynchronous generator is a transformer with a rotating secondary winding as shown in Fig \ref{fig:circuito_maquina_induccion}.  More details of this formulation can be found in \cite{machowsky}.

Let us consider the machine is integrated to an infinite bus though a transmission line as shown in Fig \ref{fig:circuito_maquina_induccion}. The losses in the line are neglected. The machine requires reactive power from the grid even if operated as generator. Usually a capacitive compensation is placed in terminals of the machine in order to minimize the voltage stability problems.  This capacitance is represented as $x_{c}$ in the figure \ref{fig:circuito_maquina_induccion}(a). A Thevening equivalent of the non-resistive components is developed in order to simplify the model as depicted in Fig \ref{fig:circuito_maquina_induccion}(b).

The equivalent Thevening impedance is given by

\begin{equation}
x_{th} = \frac{x_{m}\cdot x_{l}}{x_{m}+x_{l}} + x_{s} + x_{r}
\end{equation}

where $x_{m}$ includes the magnetizing impedance of the machine ($x_{m}'$) and the compensation device $x_{c}$

\begin{equation}
	x_{m} = \frac{x_{m}'\cdot x_{c}}{x_{m}'+x_{c}}
\end{equation}

The equivalent voltage  is 

\begin{equation}
	v_{th} = v_{B} \cdot \frac{x_{m}}{x_{m}+x_{l}}
\end{equation}

By  these considerations it is possible to determine the electrical power in the rotor  as follows:

\begin{equation}
	T_{e} = \frac{v_{th}^2\cdot r_{r}/s}{\left(r_{s} - r_{r}/s\right)^2+x_{th}^2}
\end{equation}

Using muPad we can plot this torque equation:

\begin{lstlisting}[frame=trBL]
Rs := 0.031:
Xs := 0.10:
Xm := 3.1:
Rr := 0.018:
Xr := 0.18:
Xeq := Xs + Xr:
Req := Rs + Rr:
xline := 0.08:

vth := Xm/(Xm+xline);
xth := Xm*xline/(Xm+xline) + Xs + Xr;
smax := Rr/sqrt(Rs^2+xth^2);
assume(s<smax):
Te  := vth^2*Rr/((Rs-Rr/s)^2 + xth^2)/s:
plot(Te, s=-0.4..0.4, GridVisible=TRUE):
\end{lstlisting}

\begin{figure}
	\centering
		\includegraphics[width=0.50\textwidth]{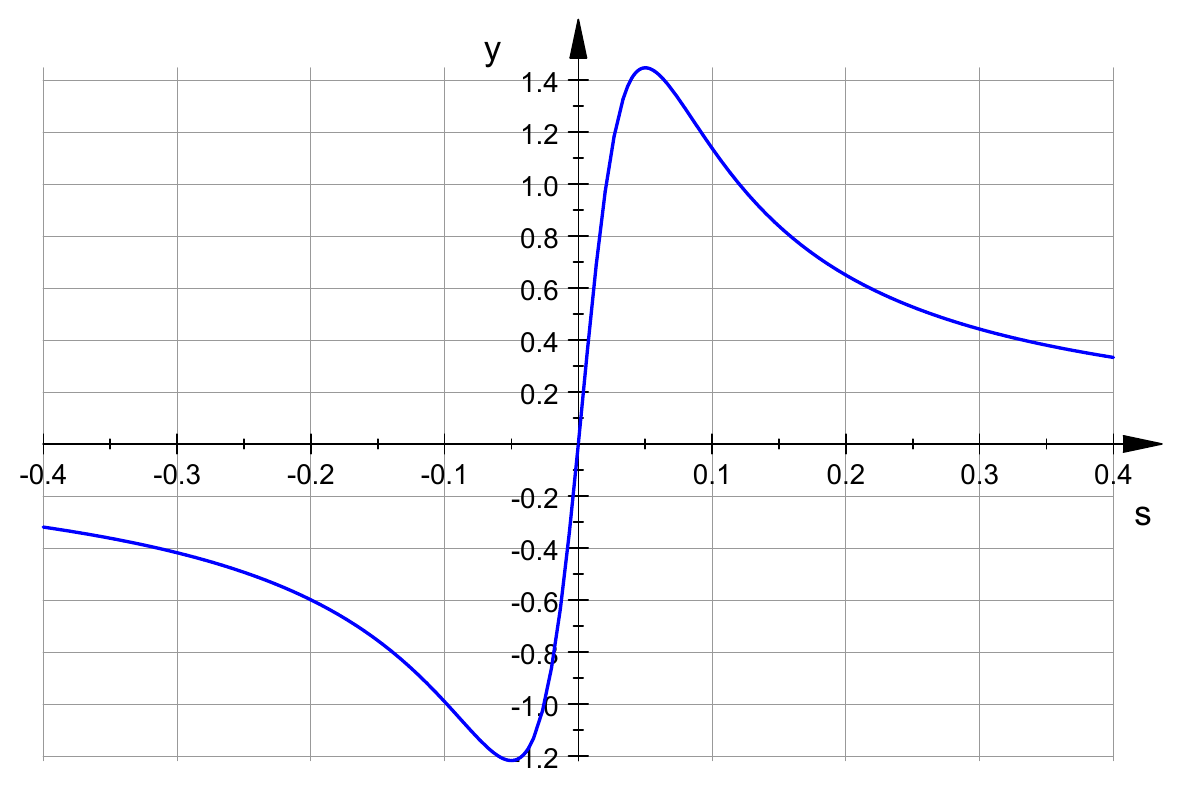}
	\caption{Electrical torque in an induction generator}
	\label{fig:Fig_Torque}
\end{figure}

Fig \ref{fig:Fig_Torque} shows the electrical torque for a particular set of parameters.  Notice the torque is positive for $s>0$ which means generation operation and negative for $s<0$ meaning motor operation.  

The maximum torque can found as the point where $\partial T_{e}/\partial s = 0$, i.e:

\begin{equation}
	s_{max} = \frac{r_{r}}{\sqrt{r_{s}^2+x_{th}^2}}
\end{equation}

At this point we are interested in the generator operation into the range $\mathcal{D}= \left\{ 0 , s_{max} \right\}$.   
\subsection{Turbine}

The power generated by a wind turbine is proportional to the cube of the wind velocity as follows:

\begin{equation}
	P = \frac{1}{2} \rho A C_{p}(\lambda)\cdot V_{w}^3
\end{equation}

where $V_{w}$ is the wind velocity, $\rho$ is the air density, $A$ the area swept by the blades, $\lambda$ the tip ratio and $C_{p}$ is the coefficient of power. The power in per unit can be expressed as follows  \cite{modeloturbina}:

\begin{equation}
	P(s) = V_{w}^3 \cdot \left(\frac{a}{\lambda} - b\right)\cdot e^{-c/\lambda} 
\end{equation}

with $w_{r(pu)} = s + 1$ because, as aforementioned $w_{s(pu)}=1$.  

The tip ratio ($\lambda$) is an important variable to take into account. It is defined as:

\begin{equation}
\lambda = \lambda_{0}\frac{V_{w}}{s+1}
\end{equation}

where $\lambda_{0}$ is the optimal tip ratio  (i.e $C_{p}=C_{p(opt)}$ when $\lambda=\lambda_0$). Notice that $P > 0$ for 

\begin{equation}
	\frac{a\cdot(s+1)}{\lambda_{0}\cdot V_{w}} - b \geq 0
\end{equation}

which fulfills in all the practical cases. In addition, the function is decreasing for $s>0$.    

Once the power is defined in terms of the slip, it is possible to define also the mechanical torque.

\begin{equation}
	T_{m} = \frac{P(s)}{s+1}
\end{equation}

\begin{figure}
	\centering
		\includegraphics[width=0.50\textwidth]{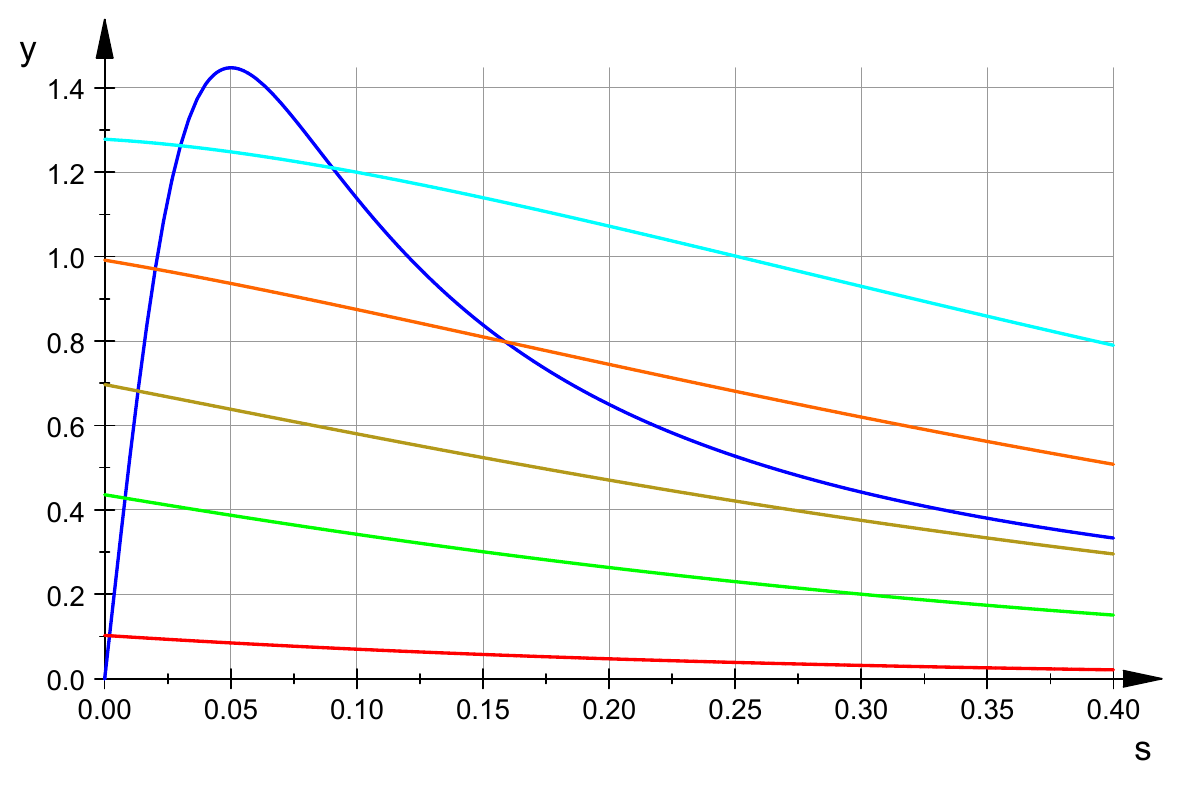}
	\caption{Electrical and mechanical torque in the turbine for different values of wind velocity and slip}
	\label{fig:Fig_TorqueMecanico}
\end{figure}

Figure \ref{fig:Fig_TorqueMecanico} depicts both, mechanical and electrical torque for different values of wind velocity and slip.  The corresponding code in muPad is shown below (notice that $T_{e}$ was already calculated).

\begin{lstlisting}[frame=trBL]
lambda_0 := 7.04:
a := 247.7079:
b := 21.6539:
c := 18.40:

Lambda :=  lambda_0*V/(1+s);
Cp     := (a/Lambda - b)*exp(-c/Lambda);
Pm     := Cp*V^3;
Tm     := Pm/(s+1):
plot(Te, Tm  $ V in [0.6,0.8,0.9,1.0,1.1], s=0..0.4,GridVisible=TRUE):
\end{lstlisting}

The resulting dynamic is depending on the slip considering $dw_{r(pu)}=ds$:

\begin{equation}
	M \frac{ds}{dt} = T_{m}(s) - T_{e}(s)
\end{equation}

where $M$ is the inertia of the turbine and generator. There are two equilibrium points (i.e points in which $T_{m}=T_{e}$ for a given wind velocity.  However, there is only one equilibrium for $s \in \mathcal{D}$ which is the interval under study.

\section{Lyapunov stability}

Let us consider a particular operative point $s_{0} \in \mathcal{D}$ and define a function $f:\mathcal{D}\rightarrow \mathbb{R}$ as:

\begin{equation}
	f(x) = \frac{1}{M}\left(T_m(x+s_0) - T_E(x+s_0)\right)
\end{equation}

Then, the dynamic system can be represented given below with the stationary point corresponding to $x=0$.

\begin{equation}
	\dot{x} = f(x)
\end{equation}

We can see the form of this function for different values of the wind velocity:

\begin{lstlisting}[frame=trBL]
M := 0.8:
xo := solve(Tm=Te,s=0..smax):
so := xo[1]:
f  := subs((Tm-Te)/M,s=x+so);
plot(f  $ V in [0.6,0.8,1.0],x=-0.05..0.05, GridVisible=TRUE): 
\end{lstlisting}

The resulting plot is depicted in Fig \ref{fig:Fig_FWind}.  Notice that the function is odd, at least in an interval $D = \left\{-s_{max},s_{max}\right\}$ which contains the $0$.
Due to this characteristic, an integral Lyapunov function $L:\mathcal{D}\rightarrow \mathbb{R}$ can be defined:

\begin{figure}[tbh]
	\centering
		\includegraphics[width=0.50\textwidth]{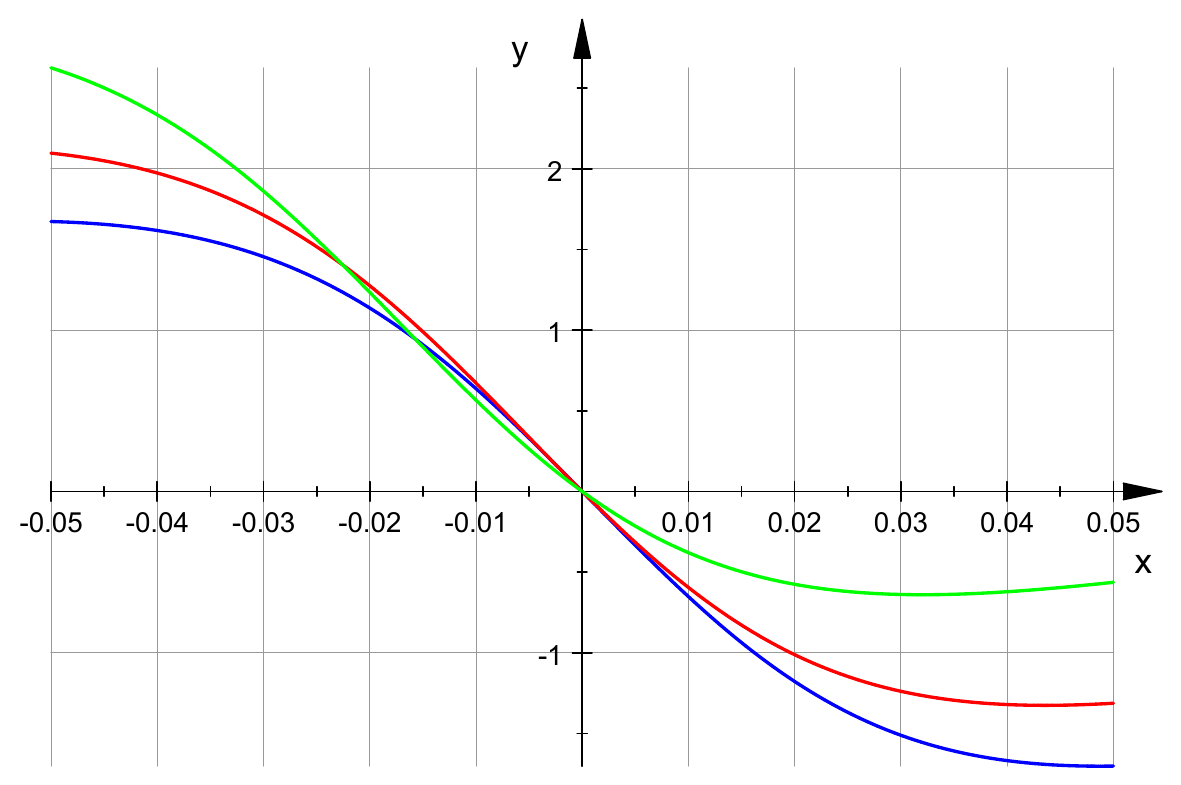}
	\caption{Function $f(x)$ for different values of wind velocity}
	\label{fig:Fig_FWind}
\end{figure}

\begin{equation}
	L(x) = -\int\limits_{0}^{x} f(x') dx'
\end{equation}

\begin{figure}[tbh]
	\centering
		\includegraphics[width=0.50\textwidth]{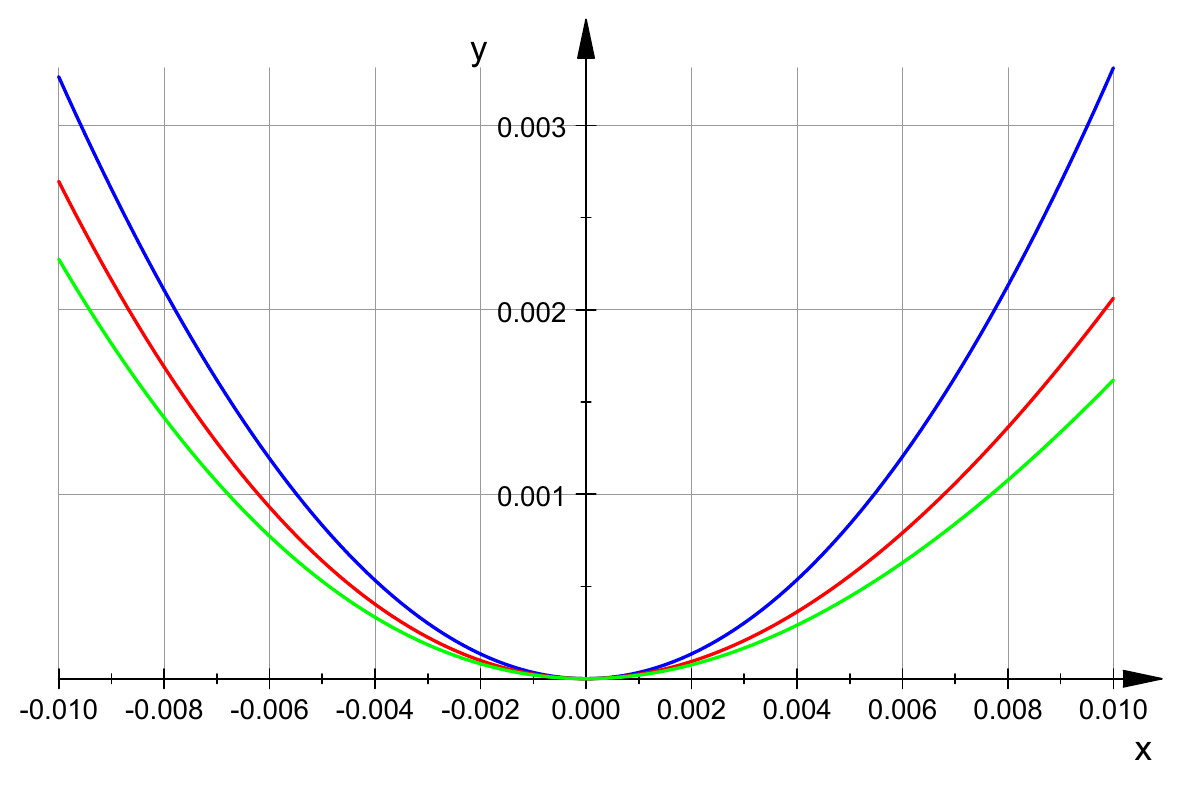}
	\caption{Proposed Lyapunov function for the wind turbine}
	\label{fig:Fig_LyapunovWind}
\end{figure}

This function is semidefinite positive $L(x)>0$ for $x \in \mathcal{D}$ and $L(0)=0$.  Therefore, it fulfills the two main conditions for being a Lyapunov function.  A plot in muPad can help to visualize these two conditions:

\begin{lstlisting}[frame=trBL]
L := -int(f,x=0..x);
plot(L  $ V in [0.6,1.0,1.05],x=-0.01..0.01, GridVisible=TRUE): 
\end{lstlisting}

Figure \ref{fig:Fig_LyapunovWind} shows the Lyapunov function for different values of wind velocity ($V_{w}$). 

Now, due to the definition as an integral function, the derivative of the function is semi-definite negative for all real values

\begin{equation}
	\frac{dL}{dt} = \frac{\partial L}{\partial x} \cdot \dot{x} = -\left( f(x)\right)^2 < 0
\end{equation}

Therefore we can conclude the system is stable.  In fact, there are only one point in which $L(0)=0$ and hence we can conclude the system is also asymptotically stable for all values of wind velocity.

\section{Region of attraction and sensitivity analysis}

\begin{figure}
	\centering
		\includegraphics[width=0.50\textwidth]{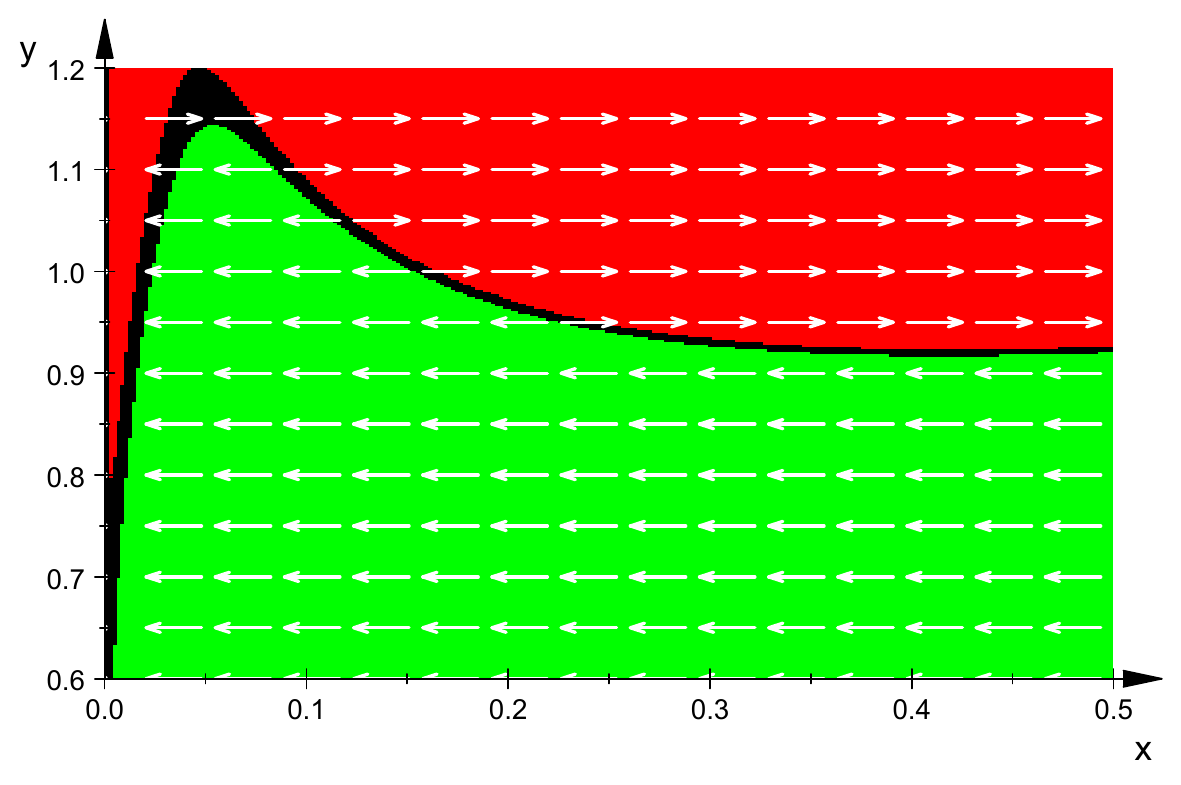}
	\caption{Region of atraction for diferent values of wind velocity}
	\label{fig:Fig_WindAtractionRegion}
\end{figure}

The attraction region can be easily determined by considering the Lyapunov function previously defined. Green region shown in Fig \ref{fig:Fig_WindAtractionRegion} shows the attraction region as function of the wind velocity.  A quiver with the direction of $f(x)$ helps to understand the phenomena.

\begin{lstlisting}[frame=trBL]
W := subs(1/M*(Tm-Te), s=x[1],V=x[2]):
SYS := (t,x) -> [W,0]:
QUIVER := plot::VectorField2d(SYS(0,x),
                              x[1] = 0..0.48, x[2] = 0.6..1.15, 
                              Mesh = [15, 12], ArrowLength=Fixed,
                              Color = RGB::White):
plot(W<0,QUIVER, x[1] = 0..0.5,x[2] = 0.6..1.2)
\end{lstlisting}

For low values of wind velocity the attractor is big enough to guarantee stability for almost any  initial condition.  However, as the wind velocity increases, the attraction region is reduced since there is other stationary point (i.e a point in which $T_{m}=T_{e}$).  This point appears in $s>s_{max}$ and is unstable since any trajectory moves away from it (see figure).   From the operation standpoint the main problem is that stable and unstable point are very close for high values of wind velocity. This can trigger instability for sudden variations of the wind velocity or the grid parameters.  	The operation of Type-A wind turbines is therefore limited to values $s<s_{max}$.  Type-A wind turbines are usually named as \textit{Fixed speed} despite not being really fixed.  This is because the narrow operational interval.

On the other hand, the compensation level modifies the stability region.  To study this phenomena we can plot the electrical torque for different compensation levels as given bellow:

\begin{lstlisting}[frame=trBL]
Xmc := Xm/(-Xm*yc+1):
vthc := Xmc/(Xmc+xline);
xthc := Xmc*xline/(Xmc+xline) + Xs + Xr;
Tec  := vthc^2*Rr/((Rs-Rr/s)^2 + xthc^2)/s:
plot(Xmc, yc = 0..3, GridVisible=TRUE):
plot(Tec $ yc = 0..3 , s=0..0.4, GridVisible=TRUE):
\end{lstlisting}

Resulting plot is depicted in Fig \ref{Fig_WindCap}.   The maximum torque increases as the compensation level increases. As consequence of that, the region of attraction is also increased.  However, in order to obtain a significant improvement it is necessary a very high compensation level.  This entails negative consequences in terms of voltage stability and costs.  

\begin{figure}
	\centering
		\includegraphics[width=0.50\textwidth]{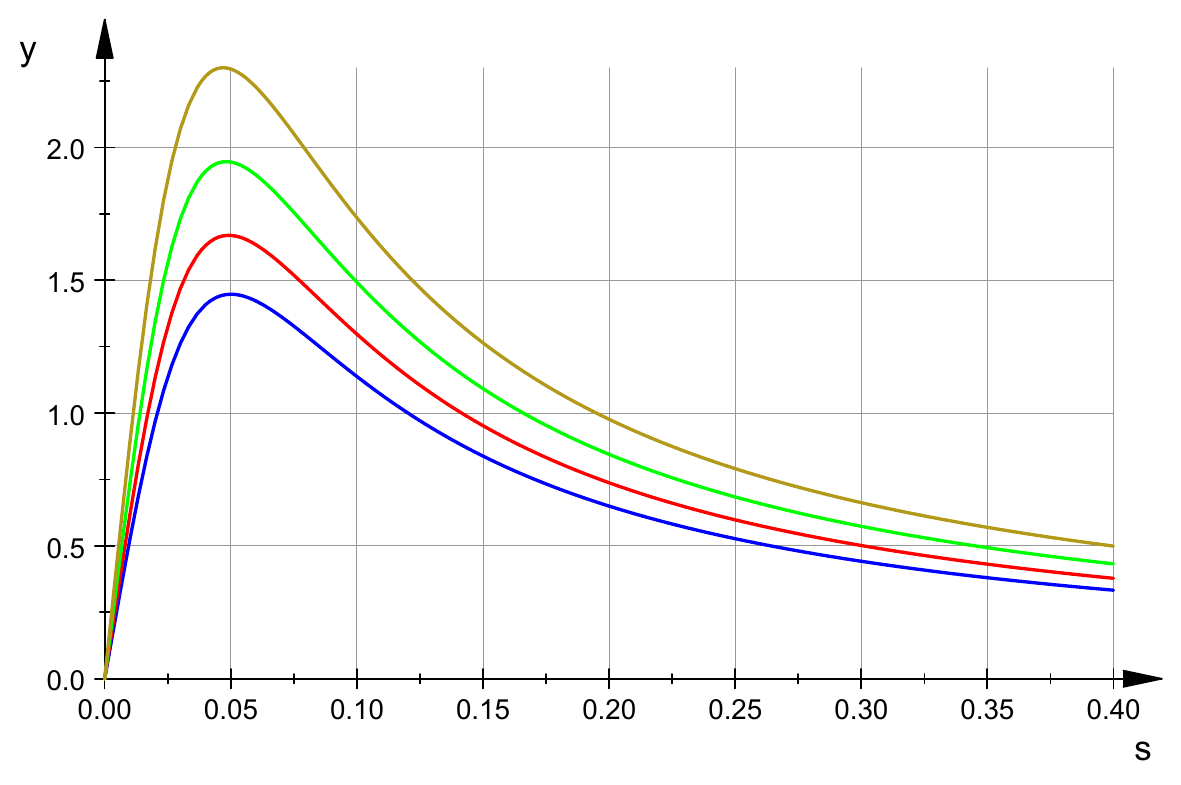}
	\caption{Electrical torque for different compensation levels}
	\label{fig:Fig_WindCap}
\end{figure}

A higher improvement on the stability is obtained by using a variable resistance in the rotor.  This variation modifies significantly the region of attraction.  Let us implement this variation in muPad for multiples of the resistance $r_r$

\begin{lstlisting}[frame=trBL]
Ten  := vth^2*r*Rr/((Rs-r*Rr/s)^2 + xth^2)/s;
plot(Ten $ r = 1..4, s=0..0.4, GridVisible=TRUE):
\end{lstlisting}

\begin{figure}
	\centering
		\includegraphics[width=0.50\textwidth]{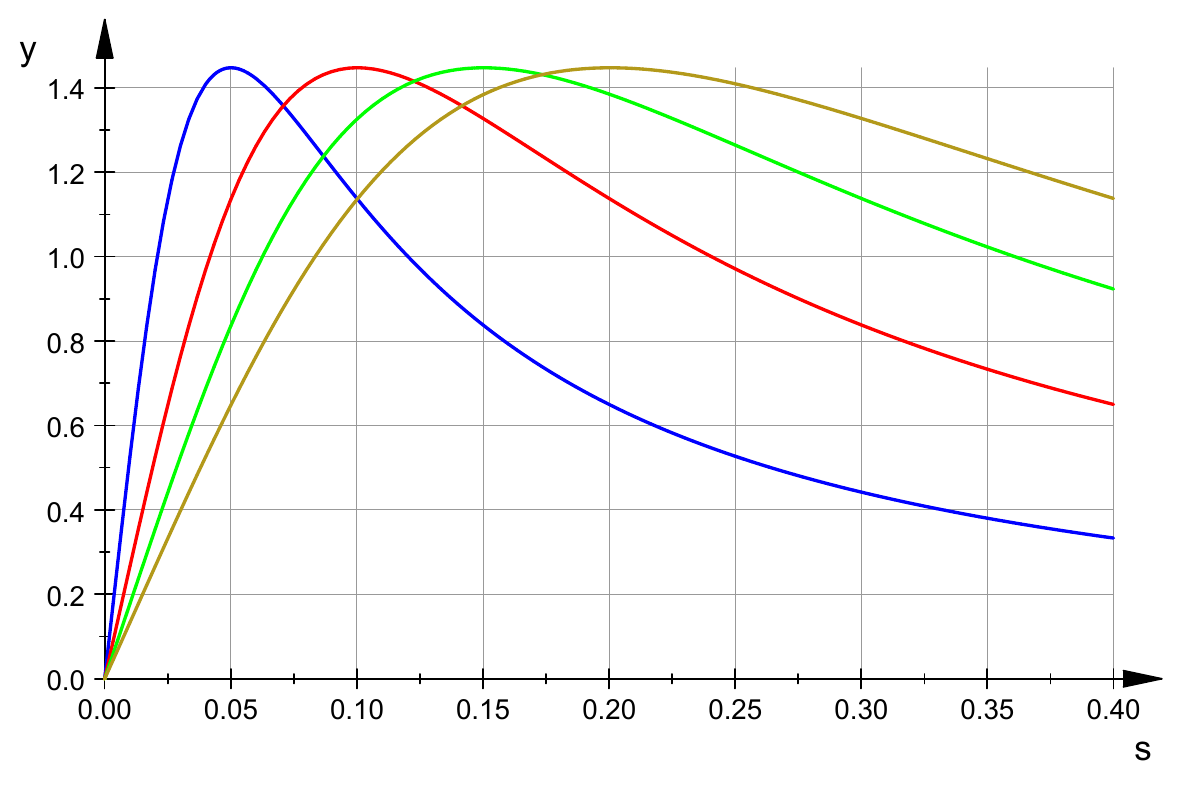}
	\caption{Electrical torque for different values of the resistance in the rotor}
	\label{fig:TipoB}
\end{figure}

Figure \ref{fig:TipoB} depicts the result of these calculations.  The region of attraction is increased even at high values of wind velocity.  However, there this imply an increasing in the losses of the machine since the rotor resistance is duplicated.      

This type of modification can be achieved by using slippery rings in order to access the rotor windings and connect a variable resistance.  This technology, usually called \textit{Type B} wind turbine, is less utilized despite its improvement in terms of stability.  The additional losses and costs are two of the main drawbacks.

\section{Review of concepts}

In this tutorial we studied the stability of a Type-A wind turbine which is a classic but current technology.  Stability was demonstrated using the Lyapunov method and considering the turbine, the induction machine and the grid.  An integral Lyapunov function was used taking into consideration the characteristic of the non-lineal equation which represents the dynamics of the system.  The impact on the stability of external capacitive compensation as well as the resistance of the rotor was considered. 
Thus, Type-B wind generators were also studied. The effect of the rotor resistance is stronger than the compensation.

\section{Further lectures}

There are a bast literature related to stability analysis of wind turbines.  In \cite{AnaliticalDFIG} a simple but efficient methodology for analysis of stability in double fed induction generator was presented.  It considers a constant mechanical torque while in reality it depends on the velocity. \cite{SatbilityPlanning}\cite{EffectParameters}\cite{Resonacia}

In \cite{LyapunovWindImportante} the Lyapunov direct method is used to study the stability of a single wind turbine. The system performance equations are linearized and expressed as simultaneous differential equations in state space.  This approach is less general than the one presented in this chapter.  However, it is very interesting in terms of concepts related to Lyapunov stability.  \cite{Damping}

\bibliographystyle{IEEEtran}
\bibliography{BibliografiaWindTA}

\end{document}